\documentclass[11pt]{article}

\usepackage[english]{babel}
\usepackage{devanagari}

\usepackage[T1]{fontenc}
\usepackage[symbol]{footmisc} 
\usepackage{mathptmx}
\usepackage{moresize}
\usepackage{pifont} 

\usepackage[affil-it]{authblk}
\usepackage{fancyhdr}
\usepackage{geometry}
\usepackage{setspace}
\usepackage{dirtytalk} 
\usepackage{changepage} 
\usepackage{filecontents} 

\MakeRobust{\say} 

\usepackage{calc} 
\usepackage{amsmath}
\usepackage{amssymb}
\usepackage[g]{esvect} 

\usepackage{enumitem} 
\setlist[itemize]{
	itemsep=-5pt,
	topsep=2pt,
	leftmargin=3.5em,
	label=\ding{226}
}

\usepackage{algpseudocode}
\usepackage[
linesnumbered,
algoruled,
vlined
]{algorithm2e}
\SetKwInput{KwInput}{Input}
\SetKwInput{KwOutput}{Output}

\usepackage{tabularray}
\usepackage{xcolor} 
\usepackage{tikz} 
\usepackage{circledsteps} 
\usepackage[percent]{overpic} 
\usepackage{adjustbox}
\usepackage[
labelfont=bf,
textfont=it,
skip=15pt
]{caption} 

\makeatletter
\g@addto@macro\@floatboxreset\centering
\makeatother 

\usepackage[
backend=bibtex,
sorting=none,
sortcites=true,
style=numeric-comp,
isbn=false,
autolang=hyphen
]{biblatex}

\AtEveryBibitem{\clearlist{language}} 
\AtEveryBibitem{
	\clearfield{urlyear}
	\clearfield{urlmonth}
}

\DefineBibliographyStrings{english}{
	url = Available,
} 

\usepackage[linktoc=all]{hyperref}
\usepackage[capitalise,nameinlink]{cleveref}
\hypersetup{
	colorlinks=true, 
	urlcolor=blue, 
	linkcolor=blue, 
	citecolor=blue, 
}

\definecolor{grey}{RGB}{128, 128, 128}

\newcommand{\brco}{\textit{breathing coefficient}\xspace}
\newcommand{\iE}{\textit{i.e.}\xspace}
\newcommand{\quadinou}{\,\,}

\newcommand{\tSup}[1]{\textsuperscript{#1}}
\newcommand{\tSub}[1]{\textsubscript{#1}}
\newcommand{\bsay}[1]{\say{\textbf{#1}}}
\newcommand{\isay}[1]{\say{\textit{#1}}}
\newcommand{\piov}[1]{\pi/#1}
\newcommand{\piovf}[1]{\frac{\pi}{#1}}

\newcommand{\xipiov}[1]{\xi_{\piov{#1}}}

\newcommand{\Rpiov}[1]{R_{\piov{#1}}}

\newcommand{\Vbody}{V_{\text{body}}}
\newcommand{\Vsolid}{V_{\text{solid}}}
\newcommand{\Vvoid}{V_{\text{void}}}
\newcommand{\DV}[1]{\Delta V_{#1}}
\newcommand{\DVbody}{\Delta \Vbody}
\newcommand{\DVsolid}{\Delta \Vsolid}
\newcommand{\DVvoid}{\Delta \Vvoid}

\newcommand{\Abody}{A_{\text{body}}}
\newcommand{\Asolid}{A_{\text{solid}}}
\newcommand{\DAbody}{\Delta \Abody}
\newcommand{\DAsolid}{\Delta \Asolid}
\newcommand{\tiAbody}{\widetilde{A}_{\text{body}}}
\newcommand{\tiAbodyo}{\widetilde{A}_{\text{body,0}}}
\newcommand{\tiAsolid}{\widetilde{A}_{\text{solid}}}
\newcommand{\tiDAbody}{\widetilde{\Delta A}_{\text{body}}}
\newcommand{\tiDAsolid}{\widetilde{\Delta A}_{\text{solid}}}

\newcommand{\udan}{\text{\dn u}}
\newcommand{\udansb}{\udan_{sb}}
\newcommand{\udanvb}{\udan_{vb}}
\newcommand{\udanbs}{\udan_{bs}}
\newcommand{\udanvs}{\udan_{vs}}
\newcommand{\udanbv}{\udan_{bv}}
\newcommand{\udansv}{\udan_{sv}}

\newcommand{\chim}{\chi_{\max}}

\newcommand{\crefalgo}[1]{\hyperref[#1]{Algorithm~\ref*{#1}}} 
\newcommand{\crefsub}[2]{\hyperref[#1]{\cref*{#1}#2}} 
\newcommand{\tmphyperref}[1]{\hyperref[label]{#1}} 

\newcommand{\sclab}[1]{#1)}
\newcommand{\sclabb}[1]{(#1)}
\newcommand{\sclbf}[1]{\textbf{#1)}}

\def\tablelinewidth{1pt}

\let\oldchi\chi
\renewcommand{\chi}{\raisebox{0.2ex}{$\oldchi$}}

\newenvironment{alignedatbracket}[1][2]{
	\left\{
	\begin{alignedat}{#1}
	}{
	\end{alignedat}
	\right.\
}

\algnewcommand{\algorithmicand}{\textnormal{\textbf{and}}\,}
\algnewcommand\AND{\algorithmicand}
\NewDocumentEnvironment{calgorithm}{O{0.5\textwidth}O{t}O{}O{}}{
	\def\tmpwidth{#1}
	
	\begin{figure*}[#2]
		\begin{center}
			\begin{minipage}{\ifx\tmpwidth\empty0.5\textwidth\else#1\fi}
				\begin{algorithm}[H]#4
					\caption{
						\rule{0pt}{10pt}\sc
						#3
					}
					\vspace{5pt}
					\DontPrintSemicolon
					\setstretch{1.2}
				}{
					\vspace{5pt}
				\end{algorithm}
			\end{minipage}
		\end{center}
	\end{figure*}
}

\newenvironment{customnomenclature}
{\begin{adjustwidth}{1em}{0pt}
		\begin{tabbing}
			\hspace{50pt}\=\kill}
		{	\end{tabbing} 
\end{adjustwidth}}

\title{Evaluation of Volume Variation Partition: \\
	the Breathing Coefficient and its Application on \\
	Uniaxial Monosized Disc Packing Swelling
	\footnote{
		\includegraphics[height=8pt,trim=0 2pt 0 0]{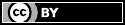} Distributed 
		under 
		a 
		\href{https://creativecommons.org/licenses/by/4.0/}{CC-BY 4.0 licence}.
	}
}

\author[]{
	Théo Boivin
	\hspace{-5pt}
	\footnote{
		Email: \href{mailto:theo.boivin@email.fr}{theo.boivin@email.fr}
	}
	\hspace{-2pt}
}

\author[]{
	\hspace{-5pt}
	Olivier Gillia
	\hspace{-5pt}
	\footnote{
		Email: \href{mailto:olivier.gillia@cea.fr}{olivier.gillia@cea.fr}
	}
}

\affil[]{Univ. Grenoble Alpes, CEA-Liten, Grenoble, F-38054, France}

\date{}

\begin{document}

\maketitle

\vspace{-1cm}

\begin{abstract}
	An analysis of the general concept of volume variation partition of a porous
	body is presented, introducing the \brco, defined as the ratio of two volume
	variations. Considering a total volume of a porous body, composed of solid
	volume and \say{void} volume, this ratio can be used to evaluate the 
	distribution of a volume variation into both others. A full description of 
	its physical interpretation is detailed, together with an uncertainty 
	analysis that specifies	precautions about its use. As an example of 
	application, a case study of 2D	monosized disc packing swelling is 
	developed. The analytical model reveals the presence of minimisation points 
	of the \brco dependent on the initial granular organisation, showing 
	possible ways to minimise the breathing of a granular material.
\end{abstract}

\vspace{0.5cm}

\hspace{8pt}
\includegraphics[width=0.87\textwidth]{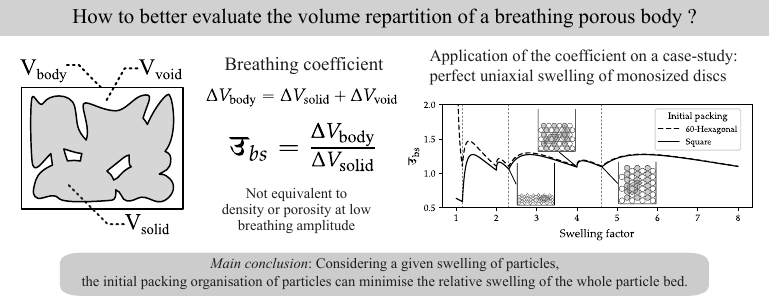}

\vspace{1cm}


\pagebreak

\section{Introduction}

The behaviour of a granular media subjected to volume variation of its particles
is nowadays considered to describe the \isay{breathing} of materials
\cite{sassine_study_2018, gillia_hydride_2021, boivin_breathing_2024}. This led
us to consider the analysis of volume variation in a more general way, not
necessarily for a granular media, but for a generic \say{porous body}. A
\isay{porous body} is defined here as a system composed of two phases: a
\say{solid} phase and a \say{void} phase (filled with gas, liquid, or void). The
term \isay{breathing} is used here to designate a succession of swelling and
shrinking of the porous body. When the need is to manage the breathing behaviour
of such a body, a major interest appears concerning the way the volume variation
of the solid phase is converted into variation of void volume or variation of
the global volume of the porous body (or similarly the way void is converted to
solid and global, or the way global is converted to void and solid). The common
quantities used to measure the relative volume fractions are the porosity
$\epsilon = \Vvoid / \Vbody$ or the solid fraction $\chi = \Vsolid / \Vbody$
(also called \say{density}), where $\Vvoid$, $\Vsolid$, and $\Vbody$ are
respectively the volumes of the void phase, the solid phase, and the porous 
body. They are really useful to describe efficiently a particular state at a 
given instant. However, how to compare two values of porosity or solid 
fractions at two different instants in time? If the porous body volume $\Vbody$ 
does not change, there is no particular issue. However, if it changes too, the 
comparison of both instants becomes unclear.

By restraining the idea of a porous body to the case of a packed particle bed,
the question of volume partition through porosity and solid fraction is well
studied. In this sense, the circle packing community has widely addressed the
static state of a circle packing \cite{cumberland_packing_1987,
bowers_introduction_2008, hifi_literature_2009}. A common methodology to study
these circle packing questions is to use a computational \say{compressor}
technique that incrementally reaches a dense state
\cite{lubachevsky_minimum_2009, specht_high_2013, markot_improved_2021,
amore_circle_2023}. Still, the attention is mostly paid to the final state but
not the intermediary behaviour. The proposed axis in this article is to extend
the notion of porosity and solid fraction to their derivative equivalents: using
volume variation instead of instantaneous volumes. As a source of inspiration,
Gomadam \& Weidner (2005) \cite{gomadam_modeling_2005} introduced the
\isay{swelling coefficient} to assess the volume variation of a porous electrode
volume and also its porosity change.

The aim of this paper is to reduce the \isay{swelling coefficient} to its
simplest form in order to generalise it as the \isay{breathing coefficient}.
\cref{sec:introduction-to-the-breathing-coefficient} introduces its definition
and different guidelines for interpretation. As an illustration,
\cref{sec:case-study-of-breathing-coefficient} develops an analytical case study
of a perfect uniaxial swelling of monosized disc packing.

\section{Introduction to the \brco}
\label{sec:introduction-to-the-breathing-coefficient}

\subsection{Definitions}

The \brco is a proposal for a basic evaluation of breathing transmission of a
porous body. It is based on the general partition of volume, divided into a
variation of solid and void volumes, where \say{void} designates the absence of
solid (filled with gas, liquid, or void). The fundamental partition is:

\begin{equation}
	\Vbody = \Vsolid + \Vvoid
\end{equation}
The partition of volume variations between two different states of the body
directly comes:
\begin{equation}
	\DVbody = \DVsolid + \DVvoid
\end{equation}

The main question approached here is to evaluate the proportion of each
contribution ($\DVsolid$ and $\DVvoid$) to the global breathing of the body
($\DVbody$), in order to assess the balance of volume variation. To this aim,
the proposed \say{\brco}, noted $\udan$ (Devan\=agar\=\i{} vowel pronounced
[u]), is defined as:
\begin{equation}
	\udan = \frac{\partial V_1}{\partial V_2}
\end{equation}
Where $V_1$ and $V_2$ can be $\Vbody$, $\Vsolid$ or $\Vvoid$. Of course, in
practice, this form is not appropriate when no analytical form is available for
$V_1$. Thus, a simplified volume variations ratio can be used as a derivative
equivalent:
\begin{equation}
	\udan = \frac{\DV{1}}{\DV{2}}
	\label{equ:udan}
\end{equation}

\begin{figure}[t!]
	\includegraphics[width=12cm]{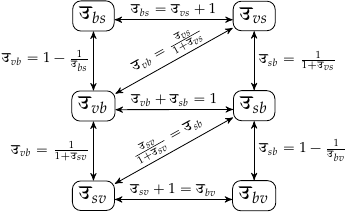}
	\caption{
		Relationships between the different definitions of the \brco.
	}
	\label{fig:relationships}
\end{figure}

More generally, six definitions exist:
\begin{equation}
	\begin{aligned}
		\udanbs = \frac{\DVbody}{\DVsolid} \quad ,
		\quad \udanvs = \frac{\DVvoid}{\DVsolid} \\[1ex]
		\udanvb = \frac{\DVvoid}{\DVbody} \quad ,
		\quad \udansb = \frac{\DVsolid}{\DVbody} \\[1ex]
		\udansv = \frac{\DVsolid}{\DVvoid} \quad ,
		\quad \udanbv = \frac{\DVbody}{\DVvoid}
	\end{aligned}
\end{equation}
Of course, these definitions remain quite equivalent, and the choice of one 
among them depends on the context and the objective of the evaluation. The 
proposedinventory is intended to be exhaustive, in order to get a complete 
toolbox. For a quick conversion, \cref{fig:relationships} details most of the 
relationships between these different definitions. We can notice that, 
mathematically, the solid and void cases are interchangeable ($v$ can be $s$ 
and reciprocally), because they remain both terms of a unique sum.

\begin{figure}[t!]
	\includegraphics[width=\textwidth]{fig_2.pdf}
	\caption{
		Different configurations and limit cases of volume variation partitions. 
		These schemes must be read relatively; the reference magnitude of 
		$\DVbody$ is just for readability (the magnitude can be close to zero). 
		Only the relative orientations matter (no negative orientation).
	}
	\label{fig:configurations}
\end{figure}

To get a more visual representation of the volume variation partition,
\cref{fig:configurations} displays a fast technique to represent the different
possible breathing configurations via a vector graph. It must be underlined that
only relative orientations matter: no negative orientation is represented. This
representation sheds light on the existence of three different configurations of
breathing which we have chosen to name: \bsay{opposite solid}, \bsay{oriented},
and \bsay{opposite void}. Furthermore, several limit cases of these
configurations exist: \bsay{void breathing}, \bsay{balanced breathing},
\bsay{solid breathing}, and \bsay{internal transfer}. In fact, the \say{internal
transfer} can be subdivided into two different cases, depending on the 
signs of $\DVsolid$ and $\DVvoid$: \bsay{solid disappearance} ($\DVsolid < 0$ 
and $\DVvoid > 0$) and \bsay{void disappearance} ($\DVsolid > 0$ and $\DVvoid < 
0$).

The \brco is theoretically a real number. Thus, it can get any value between
$-\infty$ to $+\infty$. Its interpretation is mainly based on three major
values: $0$, $1$, and $\pm\infty$. The values $0.5$ and $2$ are possibly
interesting to introduce balanced breathing, but they are not generalisable
to all \brco definitions. To help in the interpretation of this number,
\cref{tab:breathing-abacus} associates the different limit cases to the
different definitions and values.

\renewcommand{\tmphyperref}[1]{\hyperref[fig:configurations]{#1}}

\begin{table}[t!]
	\begin{tblr}{
			colspec={
				Q[c,m,0.4cm]
				Q[c,m,1cm]
				Q[c,m,1.5cm]
				Q[c,m,1.5cm]
				Q[c,m,1.5cm]
				Q[c,m,1.5cm]
				Q[c,m,1.5cm]
				Q[c,m,1.5cm]},
			row{1}={abovesep+=2pt, font=\large},
			hline{2,3}={3-Z}{\tablelinewidth},
			vline{2,3}={3-Z}{\tablelinewidth},
			cell{1}{3}={c=6}{c},
			cell{3}{1}={r=6}{c},
			stretch=1.4,
		}
		& & $\udan$ \\
		& & $-\infty$ & $0$  & $0.5$ & $1$   & $2$   & $+\infty$ \\
		\rotatebox{90}{Definition}
		& bs & VBr    & InTr & –     & SBr   & BalBr & VBr \\
		& vs & VBr    & SBr  & –     & BalBr & –     & VBr \\
		& vb & InTr   & SBr  & BalBr & VBr   & –     & InTr \\
		& sb & InTr   & VBr  & BalBr & SBr   & –     & InTr \\
		& sv & SBr    & VBr  & –     & BalBr & –     & SBr \\
		& bv & SBr    & InTr & –     & VBr   & BalBr & SBr \\
	\end{tblr}
	\vspace{-5pt}
	\caption{
		Abacus of breathing limit cases —
		\textbf{InTr}: Internal transfer,
		\textbf{VBr}: Void breathing,
		\textbf{BalBr}: Balanced breathing,
		\textbf{SBr}: Solid breathing.
	}
	\label{tab:breathing-abacus}
\end{table}

\subsection{Uncertainty analysis}

As with any other ratio number, precaution should be taken concerning the
uncertainty of the \brco. Indeed, the uncertainties due to computation and data
dumping can generate noise, and the division can enhance this noise. To
understand this issue, let us consider that the multiple mistakes due to
computation, randomness generation, and data dumping create a range of
uncertainty $\sigma \in \mathbb{R}_+$, with the same dimension as $\Delta V$.
Let us also define $\gamma_1, \gamma_2 \in [0, 1]$, two variables that are
impossible to predict due to uncertainty. Instead of \cref{equ:udan}, we obtain
in reality:
\begin{equation}
	\udan = \frac{\DV{1} \pm \gamma_1\sigma}{\DV{2} \pm \gamma_2\sigma}
\end{equation}
If the variations of volumes are equal to $\lambda_1 \sigma$ and $\lambda_2
\sigma$ where $|\lambda_1|,|\lambda_2| \gg \text{1}$, we have:
\begin{equation}
	\udan
	= \frac{\lambda_1\sigma \pm \gamma_1\sigma}{\lambda_2\sigma \pm \gamma_2\sigma}
	\approx \frac{\lambda_1\sigma}{\lambda_2\sigma}
	= \frac{\DV{1}}{\DV{2}}
\end{equation}
Here, no particular issue occurs. The uncertainty remains negligible, and the
original definition is respected. In addition, if $\DV{1}$ is now equal to
$\mu_1\sigma$ where $|\mu_1| \ll 1$, it comes:
\begin{equation}
	\udan
	= \frac{\mu_1\sigma \pm \gamma_1\sigma}{\lambda_2\sigma \pm \gamma_2\sigma}
	\approx \frac{\mu_1 \pm \gamma_1}{\lambda_2} \to 0
\end{equation}
In this case, there is still no particular issue. However, reciprocally, if
$\DV{2}$ is equal to $\mu_2\sigma$ where $|\mu_2| \ll 1$, it follows:
\begin{equation}
	\udan
	= \frac{\lambda_1\sigma \pm \gamma_1\sigma}{\mu_2\sigma \pm \gamma_2\sigma}
	\approx \frac{\lambda_1}{\mu_2 \pm \gamma_2} \to \pm\infty
\end{equation}
This time, an uncertainty appears in the final result, more precisely about the
sign. The value goes to $\infty$, but the sign is highly sensitive to the
unpredictable value $\gamma_2$. It means that when $\DV{2}$ is of the order
of magnitude of the uncertainty, it is impossible to make any conclusion about 
the real sign of $\DV{1}$ relative to $\DV{2}$ only by analysing the \brco.
Moreover, it is especially problematic when both volume variations $\DV{1}$ and
$\DV{2}$ have orders of magnitude of uncertainty, \iE are respectively equal
to $\mu_1\sigma$ and $\mu_2\sigma$:
\begin{equation}
	\udan
	= \frac{\mu_1\sigma \pm \gamma_1\sigma}{\mu_2\sigma \pm \gamma_2\sigma}
	= \frac{\mu_1 \pm \gamma_1}{\mu_2 \pm \gamma_2}
\end{equation}
In this case, this is an unsolvable limit, highly sensitive to the 
unpredictable values $\gamma_1$ and $\gamma_2$. Now, not only the sign is
concerned by uncertainty, but also the value: the whole result is
non-exploitable. As a consequence, a systematic routine to check the validity of
the \brco should be applied when using it, in order to avoid any
misinterpretation of its result. To check the value validity of $\udan$
(whatever its sign), we must have, at least:
\begin{equation}
	|\DV{1}| \gg \sigma
	\quad \text{or} \quad
	|\DV{2}| \gg \sigma
\end{equation}
To check the validity of both value and sign, we must have:
\begin{equation}
	|\DV{2}| \gg \sigma
\end{equation}

\section{Case study of the \brco:
	perfect uniaxial swelling of monosized disc packing}
\label{sec:case-study-of-breathing-coefficient}

\subsection{Problem description}

\begin{figure}[t!]
	\includegraphics[width=\textwidth]{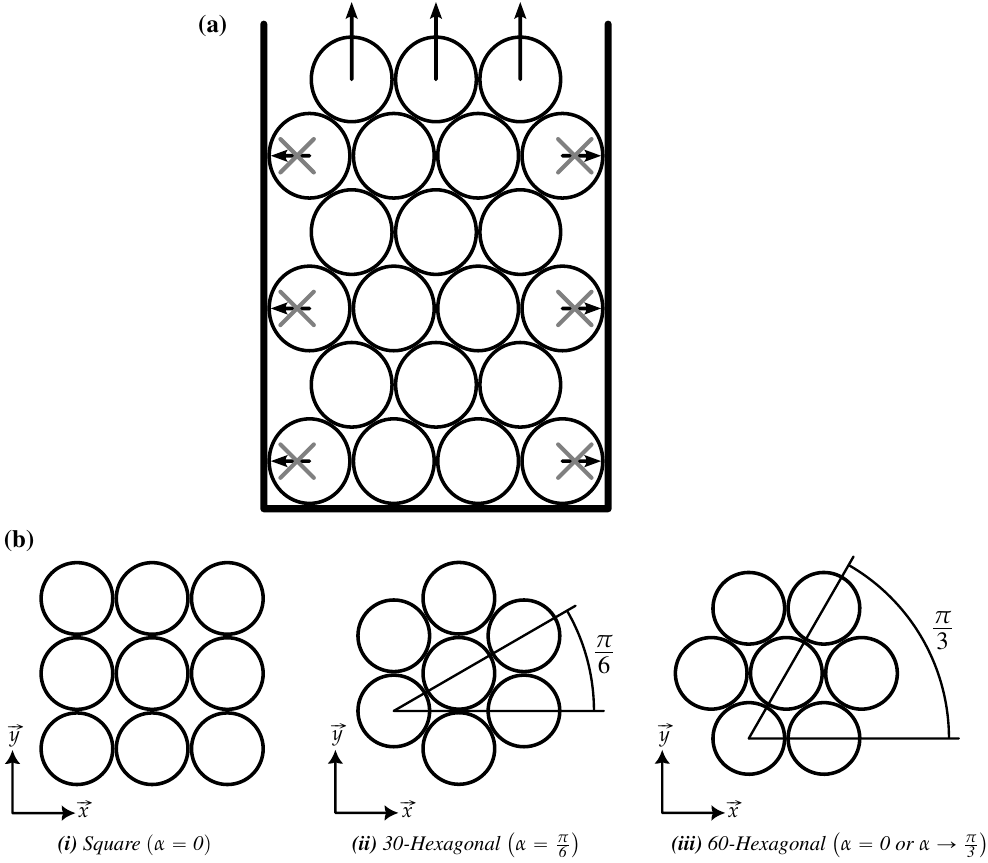}
	\caption{
		Introduction schemes of case study —
		\sclbf{a} Scheme of the problem: a container with monosized packed discs. 
		When
		the discs swell, the walls prevent the lateral movement of horizontal
		alignment going from edge to edge. The global swelling is then uniaxial.
		\sclbf{b} Three main packing types of case study.
	}
	\label{fig:container-packings}
\end{figure}

\newcommand{\tmphyperrefone}[1]	{\hyperref[fig:container-packings]{\sclab{#1}}}
\newcommand{\tmphyperreftwo}[1]	{\hyperref[fig:container-packings]{\sclabb{#1}}}
This section presents a particular application of the \brco through the
analytical modelling of a perfect uniaxial swelling of a monosized disc packing,
as schemed in \crefsub{fig:container-packings}{a}. In terms of terminology,
because the solid/void phases consideration is at the core of the \brco, the
expression \isay{disc packing} is preferred over the common expression
\say{circle packing}. However, it designates the same study area. The problem is
considered in 2D, so the \brco is based on disc area instead of volume. For this
case study, the definition of \brco is the \say{bs} one:
\begin{equation}
	\udanbs = \frac{\DAbody}{\DAsolid}
\end{equation}
The system is composed of a packing of monosized discs with initial radius
$R_0$. The swelling is defined through a swelling factor $\xi > 1$ defined as:
\begin{equation}
	\xi = \frac{R}{R_0}
\end{equation}
Where $R$ is the radius of uniformly swelled discs. The discs are perfectly
packed with a periodic scheme by neglecting the border effect. The packing is
defined through two parameters: its type and its angle $\alpha$. The type can be
\bsay{square}, \bsay{30-hexagonal} or \bsay{60-hexagonal}, as schemed in
\crefsub{fig:container-packings}{b}. The angle $\alpha$ is defined as the lowest
positive angle between horizontal and the line joining centres of two discs in
contact. Then, the square packing presents $\alpha = 0$ \tmphyperreftwo{i}
and the 30-hexagonal packing presents $\alpha = \piov{6}$ \tmphyperreftwo{ii}. 
Looking at \tmphyperrefone{iii}, the 60-hexagonal has, by
definition, $\alpha = 0$, but it also represents the limit case where $\alpha
\to \piov{3}$. Of course, the initial packing does not necessarily begin at one
of the particular packings depicted in \crefsub{fig:container-packings}{b}.
Thus, we also define an angle $\alpha_0$ that describes the initial packing
state when $\xi = 1$.

\begin{figure}[htbp!]
	\includegraphics[width=0.9\textwidth]{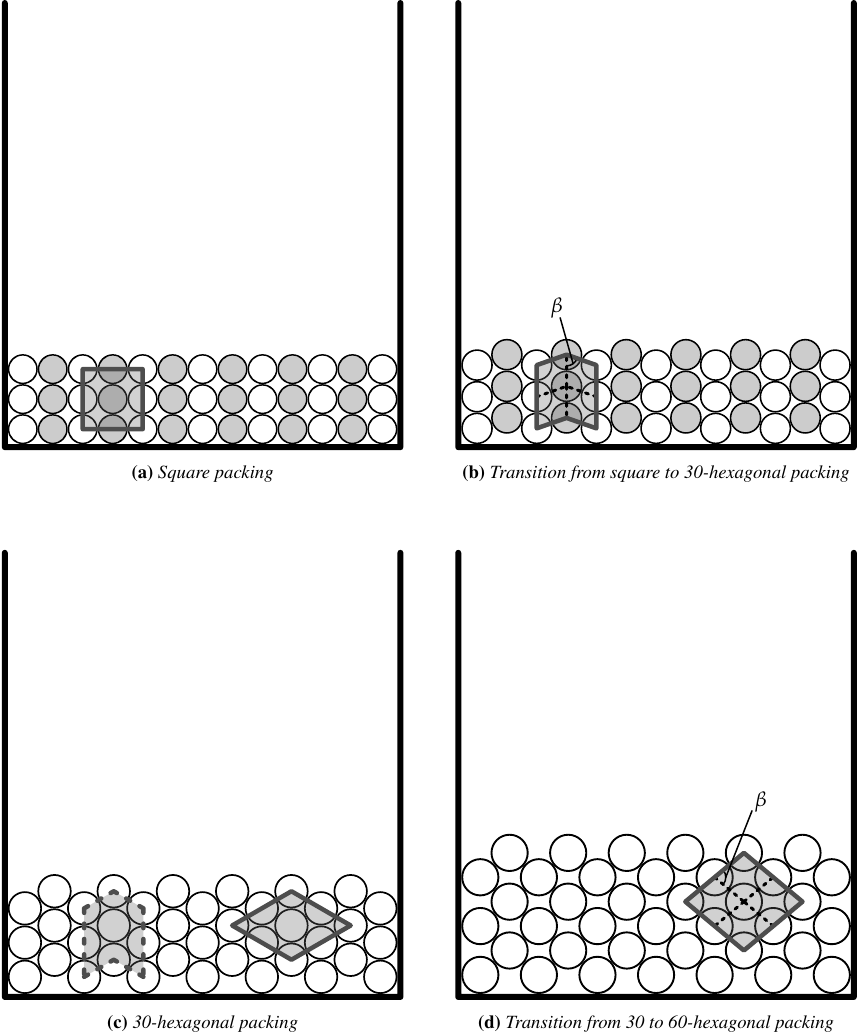}
	\caption{
		Transitions between packing types —
		Solid line polygons correspond to the representative tiles for square 
		packing in \sclbf{a} and \sclbf{b} and for 30-hexagonal packing in 
		\sclbf{c} and \sclbf{d}. Dashed polygon in \sclbf{d} corresponds to the
		representative tile of square packing after deformation due to packing 
		transition. In \sclbf{b} and \sclbf{d}, dashed lines in representative 
		tiles reveal the subtiles (four	equal parallelograms). Coloured discs 
		designate the columns that relatively move up.
	}
	\label{fig:transitions-1}
\end{figure}

\renewcommand{\tmphyperrefone}[1]	{\hyperref[fig:transitions-1]{\sclab{#1}}}
\renewcommand{\tmphyperreftwo}[1]	{\hyperref[fig:transitions-2]{\sclab{#1}}}
As introduced in \crefsub{fig:container-packings}{a}, the swelling of the
packing is uniaxial along the vertical axis ($y$ axis), which means that an
edge-to-edge horizontal alignment of discs cannot swell laterally (along $x$
axis). In this case, it necessarily leads to a granular reorganisation. This is
for this reason that the 30-hexagonal and 60-hexagonal packings, which first
seem equivalent, are in fact distinguished. The swelling consequently generates
a transition between the three main packings (square, 30-hexagonal and
60-hexagonal). The square packing being non-compact, it only occurs as an
initial packing. Then, the swelling will only be a succession of transitions
between 30-hexagonal and 60-hexagonal packings, as they are the most compact
packings when considering monosized discs \cite{chang_simple_2010}. To
understand the different transitions, \cref{fig:transitions-1,fig:transitions-2}
displays different states during swelling of discs, starting with a square
packing. Because the container constrains the packing to move vertically, the
swelling of discs generates a reorganisation. One column out of two, \iE
coloured discs in \tmphyperrefone{b}, moves up until reaching the first
30-hexagonal packing in \tmphyperrefone{c}. At this stage, all discs are free to
swell by moving up. When the swelling reaches the 60-hexagonal packing in
\tmphyperreftwo{e}, new horizontal alignments appear, which requires another
reorganisation. This reorganisation is similar to the square one, with the
difference that one row over two moves laterally over a relative distance of one
radius. It then leads to another 30-hexagonal packing, and the transition
between 30-hexagonal and 60-hexagonal begins again. For this reason, the
description of this swelling behaviour does not require a packing angle superior
to $\piov{3}$.

In a nutshell, we have $\alpha \in [0,\piov{3}[$ and $\alpha_0 = \alpha(\xi =
1).$

\makeatletter
\renewcommand{\fnum@figure}{Figure \thefigure{} (part 1)}
\makeatother

\addtocounter{figure}{-1}
\makeatletter
\renewcommand{\fnum@figure}{Figure \thefigure{} (part 2)}
\makeatother

\begin{figure}[htbp!]
	\includegraphics[width=0.9\textwidth]{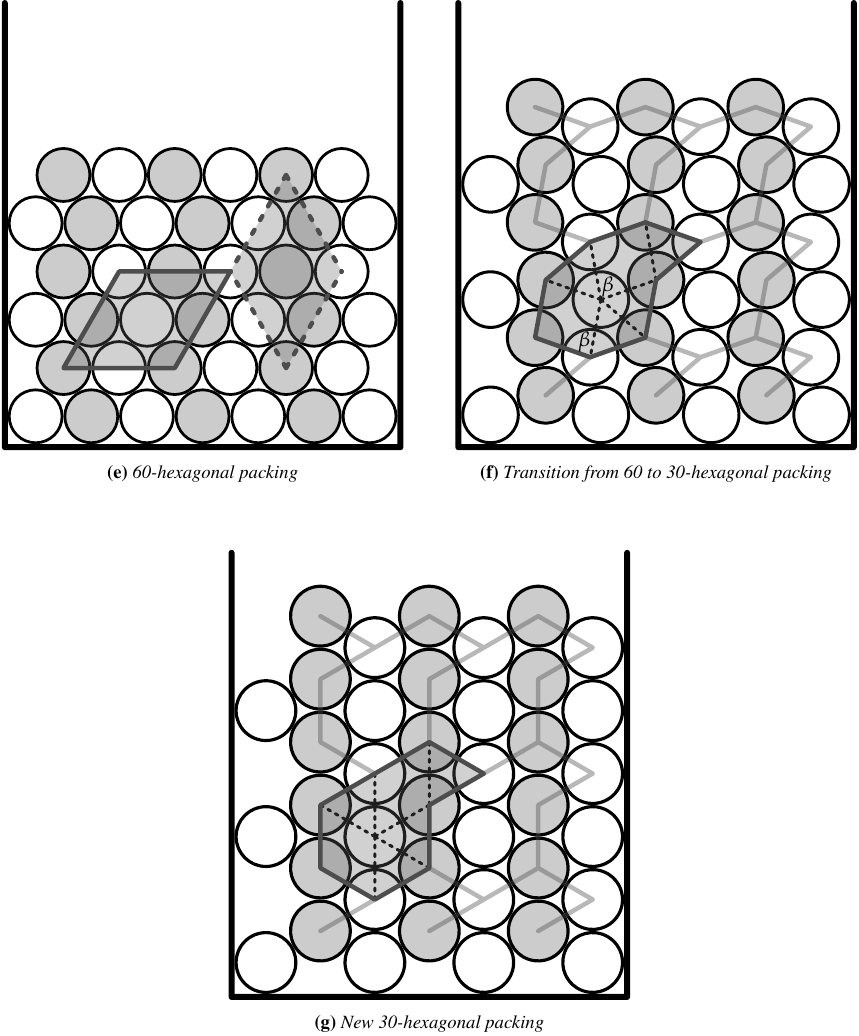}
	\caption{
		Transitions between packing types —
		Solid line polygons correspond to the representative tiles for 
		60-hexagonal packing. Dashed polygon in \sclbf{e} corresponds to the 
		representative tile of 30-hexagonal packing after deformation due to 
		packing transition. In \sclbf{f} and \sclbf{g}, dashed lines in the 
		representative tiles reveal the subtiles (two equal parallelograms and 
		four equal parallelograms). The half-transparent duplicated tiles 
		help in demonstrating the tessellation with the representative tile. 
		Coloured discs designate the columns that relatively move up.
	}
	\label{fig:transitions-2}
\end{figure}

\makeatletter
\renewcommand{\fnum@figure}{Figure \thefigure}
\makeatother

\subsection{Analytical model development}

The calculation of $\udanbs$ should be expressed only with variables $\xi$ and
$\alpha_0$, in order to nondimensionalise the resolution (independent of the 
order of magnitude of $R$). Moreover, because the problem is considered 
periodic, the calculation of areas requires the definition of representative 
tiles that should respect two conditions:
\begin{itemize}
	\item It constitutes a unique polygon sufficient to tessellate the whole 
	disc packing, only by translation, and independently of the packing angle
	$\alpha$.
	\item The areas (tile and contained discs) during transition between 
	packings are continuous.
\end{itemize}

The chosen representative tiles are the ones drawn in
\cref{fig:transitions-1,fig:transitions-2} with solid thick lines,
\tmphyperrefone{a} for square, \tmphyperrefone{c} for 30-hexagonal and
\tmphyperreftwo{e} for 60-hexagonal. The tile for square packing in
\tmphyperrefone{a} is a square that deforms at its middle. The tile for
30-hexagonal packing in \tmphyperrefone{c} is a parallelogram that keeps
constant lengths (relative to disc radii) but whose angles are varying.
Finally, the tile for 60-hexagonal packing in \tmphyperreftwo{e} requires a more
complex shape. Initially a parallelogram, it deforms differently depending on
the row.

Firstly, these tiles respect the condition of tessellation only with
translation. This is quite trivial for both first tiles, and \tmphyperreftwo{f}
and \tmphyperreftwo{g} help in demonstrating the tessellation with the third
one. Secondly, the advantage of this choice is that, whatever the value of
$\alpha$, each of these representative tiles include the area of four discs.
Whatever the angle $\alpha$, we then get:
\begin{equation}
	\Asolid = 4\pi R^2
\end{equation}
To obtain the areas of the tiles, some geometrical analytics is required by
observing the intermediary schemes, in \tmphyperrefone{a}, \tmphyperrefone{d}
and \tmphyperreftwo{f}. The first transition tile (from square to 30-hexagonal
packings) is decomposed into four parallelograms, whose symmetric areas are
directly deductible from the angle $\alpha$ (we have $\beta = \piov{2} -
\alpha$). The reasoning is equivalent for the second transition tile (from
30-hexagonal to 60-hexagonal packings), also decomposed into four
parallelograms, where $\beta = 2\alpha$. The third transition tile (from
60-hexagonal to 30-hexagonal) is decomposed into four equilateral triangles and
two parallelograms with varying angle $\beta = 2(\piov{3} - \alpha)$. Thanks to
this analysis, we obtain the respective expressions of $\Abody$ equal to the
areas of transition tiles:
\begin{equation}
	\Abody = 
	\begin{alignedatbracket}[3]
		&16 R^2 \cos(\alpha)
		\quad &&, \text{if}\quadinou \alpha \in
		\left[0, \piovf{6}\right[
		\quadinou &&\text{(square at }\alpha = 0\text{)} \\
		&4 R^2 \left[
		\sqrt{3} + 2\sin\left(\frac{2\pi}{3} - 2\alpha\right)
		\right]
		\quad &&, \text{if}\quadinou \alpha \in
		\left[0, \piovf{6}\right[
		\quadinou &&\text{(60-hexagonal at }\alpha = 0\text{)} \\
		&32 R^2 \cos(\alpha) \sin(\alpha)
		\quad &&, \text{if}\quadinou \alpha \in
		\left[\piovf{6}, \piovf{3}\right[
		\quadinou &&\\
	\end{alignedatbracket}
	\label{equ:Abody}
\end{equation}
Here again, the condition of continuity is respected. Indeed, all limits equal
$8R^2\sqrt{3}$, except for limit at $0$ for the square expression, this limit
being not a value of transition but only an initial stage never reached again
during swelling.

The last step consists in defining the relationship between the disc radii
change and the angle: $\alpha(\xi, \alpha_0)$. To do so, the swelling factor
$\xi$ is subdivided as the product of three sub-swelling factors:
\begin{equation}
	\xi = \xipiov{3} \, \xi_n \, \xi_r =
	\frac{\Rpiov{3}}{R_0} \times \frac{R_n}{\Rpiov{3}} \times
	\frac{R}{R_n}
\end{equation}
Where $\Rpiov{3}$ is the radius when $\piov{3}$ is reached for the first time
and $R_n$ the radius after $n - 1$ full transitions ($0$ to $\piov{3}$). In
other words, $\xipiov{3}$ is the swelling required to reach $\piov{3}$ for the
first time, $\xi_n$ the swelling factor required to make $n - 1$ full
transitions ($0$ to $\piov{3}$) and $\xi_r$ the remaining swelling, insufficient
to reach $\piov{3}$. Called the \bsay{transition index}, $n$ is defined as the
number of times the packing reaches the 60-hexagonal packing ($\alpha \to
\piov{3}$). We can notice that, starting from a horizontal alignment of discs
($\alpha = 0$, square or 60-hexagonal), the distance between both centres does
not change (for each line, relatively, one disc out of two moves up, but not
laterally). The only difference between square and 60-hexagonal packings is
that, in the second case, one line out of two laterally moves during the
transition to the 30-hexagonal packing. Then, the initial horizontal distance
between centres $d_h$ remains constant until $\alpha$ reaches $\piov{3}$:
\begin{equation}
	d_h = 2 \xipiov{3} R_0 \cos\left(\piovf{3}\right) = 2 R_0 \cos(\alpha_0)
	\implies \xipiov{3} = 2\cos(\alpha_0)
\end{equation}
Where, with $\alpha_0 = 0$, it comes $\xipiov{3} = 2$. By consequence, the
transition of $\alpha$ from $0$ to $\piov{3}$ (whatever it starts from a square
or a 60-hexagonal packing) requires the disc radius to be doubled. This
conclusion implies that $\xi_n$ is of the form:
\begin{equation}
	\xi_n = 2^{n(\xi,\alpha_0)-1}
\end{equation}
Where $n$ depends on $\xi$ and can be determined with the \crefalgo{alg:n}. The
remaining swelling factor is directly defined by the final packing angle
$\alpha$:
\begin{equation}
	\xi_r = \frac{1}{\cos(\alpha)} = \frac{\xi}{\xipiov{3}\xi_n}
\end{equation}
It comes the final expression of packing angle $\alpha$:
\begin{equation}
	\alpha = \arccos\left(\frac{2^{n(\xi,\alpha_0)}\cos(\alpha_0)}{\xi}\right)
	\label{equ:alpha}
\end{equation}
This leads to the final form of the \brco for this case study. Because of the
ratio, it can be simplified with nondimensionalised areas (independency from the
radius value), by dividing by $4R_0^2$:
\begin{equation}
	\udanbs = \frac{\tiDAbody}{\tiDAsolid}
\end{equation}
Where:
\begin{equation}
	\begin{alignedatbracket}
		\tiDAbody &= \tiAbody - \tiAbodyo \\
		\tiDAsolid &= \pi(\xi^2 - 1)
	\end{alignedatbracket}
	\label{equ:nondim-Deltas}
\end{equation}
The expression of $\tiDAbody$ is based on \cref{equ:Abody}, depending on the
packing angle $\alpha$, through the \crefalgo{alg:DAbody}. Finally, an
interesting common parameter is the solid fraction, directly defined based on
dimensionless areas:
\begin{equation}
	\chi = \frac{\tiAsolid}{\tiAbody}
\end{equation}

\begin{calgorithm}[0.5\textwidth][t][$n$ computation][\label{alg:n}]
	\KwInput{
		$\xi$,
		$\alpha_0$
	}
	$n = 0$\;
	$\xi' = \frac{\xi}{\cos(\alpha_0)}$\quadinou \# Reset angle to zero \;	
	\While{$\xi' > 2$}{
		$n \gets n + 1$, $\xi' \gets \frac{\xi'}{2}$\;
	}
	\vspace{2pt}
	\KwOutput{
		$n$
	}
\end{calgorithm}
\vspace*{-5pt}

\begin{calgorithm}[0.6\textwidth][t!][$\tiDAbody$ computation][\label{alg:DAbody}]
	\KwInput{
		$\xi$,
		$\alpha_0$
	}
	$\alpha =$\quadinou\cref{equ:alpha}\;
	\If{$\alpha < \piovf{6}$}{
		\If{Square initial packing \AND $n = 0$}{
			$\tiAbody = 4\xi^2\cos(\alpha)$
		}
		\Else{
			$\tiAbody = \xi^2 \left[
			\sqrt{3} + 2\sin\left(\frac{2\pi}{3} - 2\alpha\right)
			\right]$
		}
	}
	\Else{
		$\tiAbody = 8\xi^2\cos(\alpha)\sin(\alpha)$
	}
	\If{$\alpha_0 < \piovf{6}$}{
		\If{Square initial packing}{
			$\tiAbodyo = 4\cos(\alpha_0)$
		}
		\If{60-Hexagonal initial packing}{
			$\tiAbodyo = 
			\sqrt{3} + 2\sin\left(\frac{2\pi}{3} - 2\alpha_0\right)$
		}
	}
	\label{line:test}
	\Else{
		$\tiAbody = 8\cos(\alpha_0)\sin(\alpha_0)$
	}
	$\tiDAbody = \tiAbody - \tiAbodyo$\;
	\KwOutput{
		$\tiDAbody$
	}
\end{calgorithm}
\vspace*{-5pt}

\subsection{Results discussion}

\begin{figure}[t!]
	\includegraphics[width=0.8\textwidth]{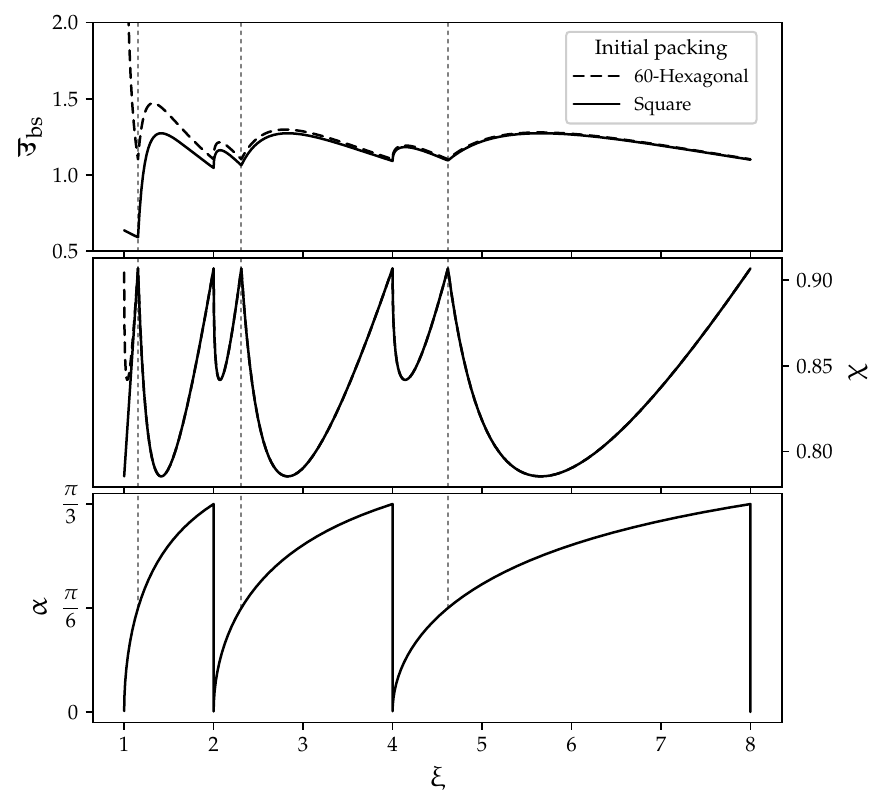}
	\caption{
		Breathing coefficient $\udanbs$, solid fraction $\chi$ and packing angle
		$\alpha$ as a function of swelling factor $\xi$ and initial packing 
		(square in solid line, 60-hexagonal in dashed line), with $\alpha_0 = 0$. 
		Vertical dotted lines corresponds to 30-hexagonal packing ($\alpha = 
		\piov{6}$).
	}
	\label{fig:curves:VSxi-dualpacking}
\end{figure}

As a particular case of the model, \cref{fig:curves:VSxi-dualpacking} displays
the evolution of the \brco $\udanbs$, the solid fraction $\chi$ and the packing
angle $\alpha$ as a function of the swelling factor $\xi$. It is more especially
restrained to the case with a null initial packing angle $\alpha_0$, for both
square and 60-hexagonal initial packings. An overview of the graph indicates an
order of magnitude of $\udanbs$ between $1$ and $1.5$, which is close to a
\say{solid breathing} (see \cref{tab:breathing-abacus} for interpretation). In
other words, most of the swelling of the body is induced by the disc swelling,
and the generation of void area is relatively low. However, a singular behaviour
appears at the beginning of the 60-hexagonal initial packing case (dashed line
in \cref{fig:curves:VSxi-dualpacking}), where $\udanbs$ tends towards infinity
for low values of swelling factor $\xi$. In this case, the swelling is mostly
due to void area variation. This is coherent with the fact that the 60-hexagonal
packing is the most compact packing for monosized discs, so at this state, the
slightest variation of disc area necessarily leads to a high variation of void.
On the opposite, the square initial packing displays a low \brco at a low
amount of swelling, presenting a limit for $\xi \to 1$ equal to $2/\pi \approx
0.634$. At this stage, the swelling is at the transition between \say{internal
transfer} and \say{solid breathing} (see \cref{tab:breathing-abacus} for
interpretation): the area variation of discs is compensated by an opposite area
variation of void. These particular behaviours changes when the swelling is
enough to reach the 30-hexagonal packing ($\alpha = \piov{6}$). For both square
and 60-hexagonal packings, the \brco suddenly increases with swelling factor
$\xi$, which is coherent, once again considering the compact characteristic of
the 30-hexagonal packing. A first transition finishes when the swelling factor
$\xi$ reaches $2$. From this point, both square and 60-hexagonal initial
packings tend to a similar behaviour, revealing \brco local minimisation points
for every $\xi$ value of $2^n$ (at each 60-hexagonal packing) and
$2^{n+1}\sqrt{3}/3$ (at each 30-hexagonal packing). In fact, these local
minimisation points can be quite interesting for any process that uses granular
breathing, where the estimation of volume variation can be important for
container design. For example, by taking the 60-hexagonal initial packing case
(dashed line in \cref{fig:curves:VSxi-dualpacking}), for a swelling factor
comprised between $2$ and $4$, the minimum value of $\udanbs$ is around $1.103$
(at $\xi = 4$) and the maximum is around $1.298$ (at $\xi = 2.802$). The
respective equivalent values by taking the definition \say{vs} of \brco are
$0.103$ and $0.298$. In other words, at $\xi = 4$, the swelling of discs is
higher than at $\xi = 2.802$, but a lower proportion of void area is generated
($10.3$\% against $29.8$\%). Consequently, an optimisation step seems necessary
for any process where, for example, the density of the granular material is
important. More precisely, looking at the solid fraction curve $\chi$, the
respective values at $\xi = 4$ and $\xi = 2.802$ are around $0.907$ and $0.785$.
The case at $\xi = 4$ is the hexagonal configuration proved to be the densest
for a monosized disc packing, with $\chim = \piov{\sqrt{12}}$,
\cite{chang_simple_2010}.

\begin{figure}[t!]
	\includegraphics[width=0.8\textwidth]{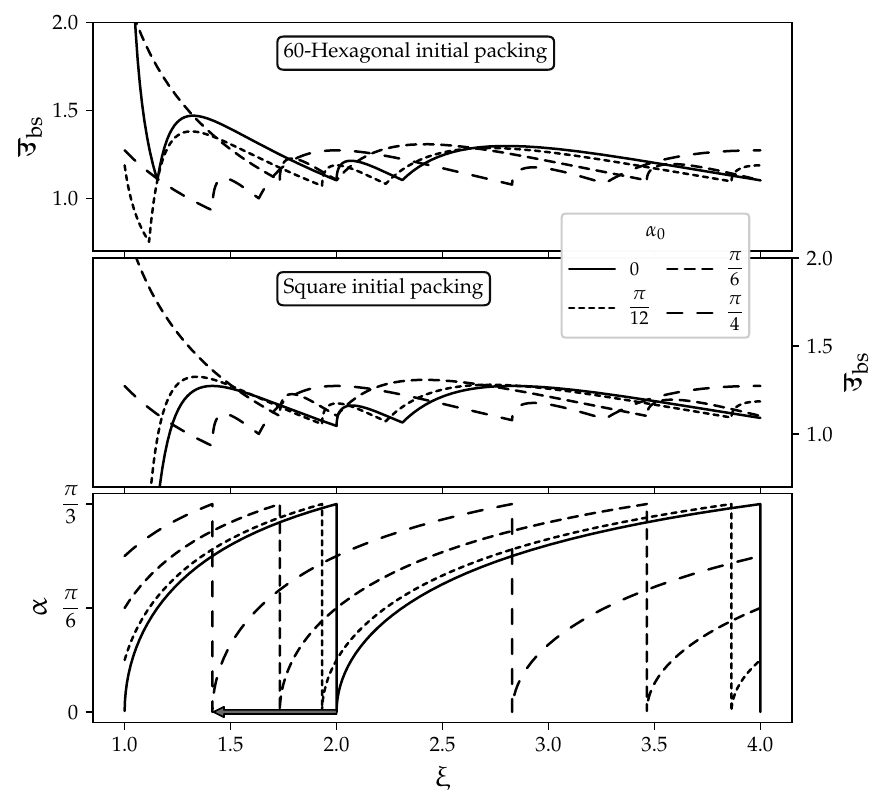}
	\caption{
		Breathing coefficient $\udanbs$ and packing angle $\alpha$ as a function 
		of swelling factor $\xi$ and initial angle $\alpha_0$. For breathing 
		coefficient, top curve is 60-hexagonal initial packing, bottom curve is 
		square initial packing. The horizontal arrow draws the phase shift due to 
		variation of initial angle $\alpha_0$.
	}
	\label{fig:curves:VSxi-alpha0}
\end{figure}

The initial angle $\alpha_0$ also has an important influence on the distribution
of local minimisation points of \brco. In this sense,
\cref{fig:curves:VSxi-alpha0} displays the \brco $\udanbs$ and the packing angle
$\alpha$ as a function of the swelling factor $\xi$, and this for several values
of initial angle $\alpha_0$. Both square and 60-hexagonal initial packings are
displayed for the \brco. Expectedly, the difference between square and
60-hexagonal initial packings only occurs for $\alpha_0 < \piov{6}$. Indeed,
above $\piov{6}$, even in the square initial packing case, the packing is
already hexagonal. As observed above, the difference is more pronounced for a
swelling factor $\xi$ comprised between $1$ and $2$, where the square initial
packing generates lower values of \brco, especially for $\alpha < \piov{6}$.
Generally speaking, for both square and 60-hexagonal initial packings, the
increase of $\alpha_0$ implies a phase shift to the left, as drawn by the arrow
on the $\alpha$ curve. Consequently, the formation of 30-hexagonal or
60-hexagonal packings does not occur at the same amount of swelling. More
precisely, these points are more generally met at every $\xi$ value of
$2^n\cos(\alpha_0)$ (at each 60-hexagonal packing) and
$2^{n+1}\cos(\alpha_0)\sqrt{3}/3$ (at each 30-hexagonal packing). This
formulation includes the particular case addressed above, for $\alpha_0 = 0$.
This observation is interesting for a practical application where the maximum
swelling factor is imposed. In this case, a judicious choice of initial granular
organisation (through $\alpha_0$) would imply a minimal \brco at the highest
amount of swelling.

As we can observe, the model is not linear with swelling factor $\xi$, but
logarithmic in base $2$ (new $\alpha$ transition every $2^n$). The local
minimisation points of the \brco of these transitions revealed to be at
hexagonal packings ($\alpha = 0$ and $\alpha = \piov{6}$). More precisely, the
expressions of the \brco at these points are, respectively:

\begin{equation}
	\begin{alignedatbracket}[2]
		\udanbs(\alpha = 0) &= \frac
		{2^{2n+1}\sqrt{3}\cos^2(\alpha_0) - \tiAbodyo}
		{\pi(2^{2n}\cos^2(\alpha_0) - 1)} \\[1ex]
		\udanbs\left(\alpha = \piovf{6}\right) &= \frac
		{2^{2n+3}\sqrt{3}\cos^2(\alpha_0) - 3\tiAbodyo}
		{\pi(2^{2n+2}\cos^2(\alpha_0) - 3)} \\
	\end{alignedatbracket}
\end{equation}
From these expressions, a particular limit value emerges at high transition
index $n$:
\begin{equation}
	\lim_{n\to\infty}\udanbs(\alpha = 0)
	= \lim_{n\to\infty}\udanbs\left(\alpha = \piovf{6}\right)
	= \frac{\sqrt{12}}{\pi}
	\approx 1.102\,657\,791
\end{equation}

\begin{figure}[t!]
	\includegraphics[width=0.8\textwidth]{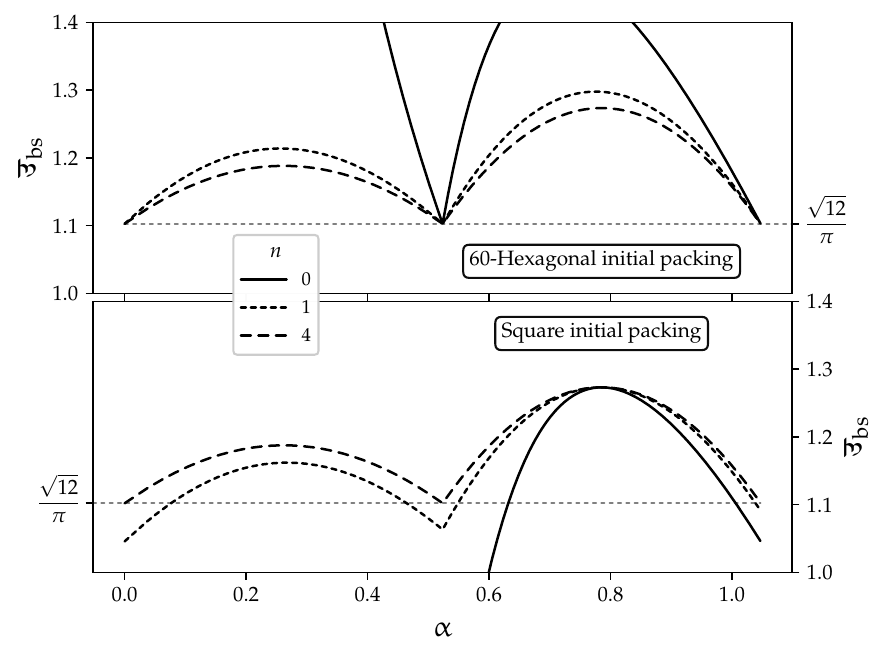}
	\caption{
		Breathing coefficient $\udanbs$ as a function of packing angle $\alpha$ 
		and initial packing. The horizontal dashed lines represent the limit 
		minimisation value of $\sqrt{12}/\pi$. The initial packing angle 
		$\alpha_0$ is null.
	}
	\label{fig:curves:VSalpha}
\end{figure}

Considering the high amount of swelling of monosized discs supposed to
reorganise in hexagonal packings, this value can be seen as the theoretical
local minimisation value of the \brco. It must be noticed that significantly
lower values of the \brco can be obtained but only at low values of swelling ($n
< 4$). However, for higher swelling ($n \ge 4$), the minimisation to
$\sqrt{12}/\pi$ becomes particularly valid. In this optic,
\cref{fig:curves:VSalpha} shows the \brco as a function of packing angle
$\alpha$ for both square and 60-hexagonal initial packings, with a null initial
packing angle $\alpha_0$, and this for different transition indexes $n$. A
difference appears between square and 60-hexagonal initial packings, where in
the square case the limit to the minimisation value $\sqrt{12}/\pi$ is made with
inferior values, whereas in the 60-hexagonal case, it is made with superior
values. In \cref{fig:curves:VSalpha}, the curves for $n > 4$ (not displayed
here) are quite merged with the curve at $n = 4$, which confirms the transition
index limit stated above. After $n = 4$, the distinction between square and
60-hexagonal initial packings becomes negligible.

In fact, this limit value of $\sqrt{12}/\pi$ is not a surprise at all, because
it corresponds to the reverse of the maximum theoretical solid fraction:
\begin{equation}
	\lim_{n\to\infty}\udanbs(\alpha = 0) = \frac{1}{\chim}
\end{equation}
This observation reveals an almost equivalence between porosity, solid fraction,
and \brco. When looking at \cref{equ:nondim-Deltas}, it is obvious that:
\begin{equation}
	\lim_{\xi\to\infty}\udanbs = \frac{1}{\chi}
\end{equation}
And also, by taking the \say{vs} definition, the equivalence is made with
porosity:
\begin{equation}
	\lim_{\xi\to\infty}\udanvs = \frac{1}{\epsilon}
\end{equation}

In other words, when the initial areas $A_0$ (or, more generally speaking,
initial volumes $V_0$) become negligible compared to their respective
variations, the \brco becomes useless, bringing nothing additional compared to
the traditional porosity and solid fraction. A first impression of equivalency
between porosity or solid fraction compared to \brco may come to mind. In
fact, this is only true when $\DV{1} \gg V_{1,0}$ and $\DV{2} \gg V_{2,0}$,
where $V_{1,0}$ and $V_{2,0}$ designates the initial volumes. In practice, it is
quite expected that the volume variations, if they are not at the same order of
magnitude, will be much lower than the absolute volumes (for example, thinking
about a thermal dilation of a particle bed or a gas storage within a hybrid
porous material). In this case, the way the volume variations between body,
solid, and void are distributed remains hard to apprehend. It is for this reason
that the \brco invites thinking in terms of volume variations instead of
absolute volumes.

\section{Conclusions}

The aim of this paper was to introduce an analytical description of the \brco,
noted $\udan$ and defined as $\DV{1} / \DV{2}$, where $\DV{1}$ and $\DV{2}$ are
two volume variations. It is mostly useful for a porous body whose total volume
changes due to a volume variation of solid and another of \say{void} (filled
with gas, liquid, or void). Based on a basic volume variation partition, this
parameter presents six different definitions, whose choice depends on the
applied problem. The interest of this number is to evaluate the conversion of a
volume variation into both others. Together with the description of an
interpretation guide, an uncertainty analysis of this parameter was made.
Considering a given uncertainty $\sigma$, we suggest a systematic routine in
order to check the validity of the \brco, with $|\DV{1}| \gg \sigma$ or
$|\DV{2}| \gg \sigma$ to validate the value, and $|\DV{2}| \gg \sigma$ to
validate both value and sign.

As an example of application, a case study of perfect uniaxial swelling of
monosized disc packing was carried out, using a 2D analytical model. In this
case study, the discs swelling generates a global swelling of the discs bed,
where periodic packing types are met. Focusing on the definition $\udanbs =
\DVbody / \DVsolid$, the analysis revealed a periodic evolution with the
swelling factor, giving birth to regular local minimisation points of the
\brco (when the disc packing becomes hexagonal). This periodicity occurs
every $2^n\cos(\alpha_0)$ and $2^{n+1}\cos(\alpha_0)\sqrt{3}/3$, by considering
the value of the swelling factor, where $\alpha_0$ is the initial packing angle.
The position of the local minimisation points of the \brco depends on
$\alpha_0$. This observation constitutes an interesting lead for the design of
porous material that would self-absorb the swelling of its solid phase. At high
swelling factor ($\ge 2^4$), the limit of minimal local point for the \brco is
found to be $\sqrt{12}/\pi$, which correspond to the reverse of the maximal
solid fraction for a monosized disc packing.

\section*{Nomenclature}

\subsection{Main symbols}
\begin{customnomenclature}
	$\udan$ \> Breathing coefficient, – \\
	$A$ \> Area, m\tSup{3} \\
	$d$ \> Distance, m \\
	$n$ \> Transition index, – \\
	$R$ \> Radius, m \\
	$V$ \> Volume, m\tSup{3} \\
	$\alpha$ \> Packing angle, – \\
	$\beta$ \> Intermediary angle, – \\
	$\gamma$ \> Undefined variable $\gamma \in [0, 1]$, – \\
	$\lambda$ \> Noticeable factor $|\lambda| \gg 1$, – \\
	$\mu$ \> Negligible factor $|\mu| \ll 1$, – \\
	$\xi$ \> Swelling factor, – \\
	$\epsilon$ \> Porosity, – \\
	$\sigma$ \> Volume uncertainty, m\tSup{3} \\
	$\chi$ \> Solid fraction, – \\
\end{customnomenclature}

\pagebreak

\subsection{Subscripts}
\begin{customnomenclature}
	0 \> Initial value (at $\xi = 1$) \\
	1 \> Numerator term \\
	2 \> Denominator term \\
	body \> Porous body \\
	$bs$ \> Body over solid \\
	$bv$ \> Body over void \\
	$h$ \> Horizontal \\
	max \> Maximum value \\
	$n$ \> Required to reach $n - 1$ full transitions \\
	solid \> Solid phase \\
	void \> \say{Void} phase (with gas, liquid, or void) \\
	$vb$ \> Void over body \\
	$vs$ \> Void over solid \\
	$r$ \> Remaining swelling \\
	$sb$ \> Solid over body \\
	$sv$ \> Solid over void \\
	$\piov{3}$ \> Required to reach $\alpha = \piov{3}$ for the first time \\
\end{customnomenclature}

\subsection{Operations}
\begin{customnomenclature}
	$\Delta y$ \> Variation, [y] \\
	$y'$ \> Secondary value, [y] \\
	$|y|$ \> Absolute value, [y] \\
	$\widetilde{y}$ \> Nondimensionalised equivalent, – \\
	$\partial y_1 / \partial y_2$ \> Partial derivative, [y\tSub{1}/y\tSub{2}] \\
\end{customnomenclature}

\section*{Acknowledgments}

The authors acknowledge the collaboration of the team from Umicore company
composed of Duancheng Ma, Michal Tulodziecki and Jacob Locke. They also
gratefully thank Benoît Mathieu and Willy Porcher for their rich involvement in
the thesis project.

\section*{Credit authorship contribution statement}

T.B. designed the figures, drafted and wrote the manuscript, O.G. supervised the
work, T.B. and O.G. developed the analytical model, discussed the results and
reviewed the manuscript.

\section*{Declaration of generative AI and AI-assisted technologies in the 
writing process}

During the preparation of this work, T.B. occasionally used
\isay{\href{https://www.deepl.com}{DeepL}} as a French-to-English translator for
unique words or short expressions (up to five words). T.B. also used
\isay{\href{https://languagetool.org}{LanguageTool}} for grammar / spelling
check. After using this tool, the authors reviewed and edited the content as
needed and take full responsibility for the content of the published article.

\section*{Declaration of competing interest}

The authors declare no conflicts of interest regarding this manuscript.

\section*{Data availability}

The data that support the findings of this study are openly available at the 
Gitlab repository \cite{gitlab_breathing}.

\printbibliography

@THESIS{sassine_study_2018,
	AUTHOR = {Sassine, Nahia},
	INSTITUTION = {Université Grenoble Alpes},
	DATE = {2018-11},
	LANGID = {english},
	TITLE = {Study of the Thermo-Mechanical Behavior of Granular Media and 
	Interactions Medium-Tank},
	TYPE = {phdthesis},
	URLDATE = {2024-10-14},
	URL = {https://theses.hal.science/tel-02043455},
}

@ARTICLE{gillia_hydride_2021,
	AUTHOR = {Gillia, O.},
	DATE = {2021},
	DOI = {10.1016/j.ijhydene.2021.07.082},
	ISSN = {0360-3199},
	JOURNALTITLE = {International Journal of Hydrogen Energy},
	LANGID = {english},
	PAGES = {35594-35640},
	SHORTTITLE = {Hydride Breathing and Its Consequence on Stresses Applied to 
	Containers},
	TITLE = {Hydride Breathing and Its Consequence on Stresses Applied to 
	Containers: {{A}} Review},
	URLDATE = {2021-09-27},
}

@ARTICLE{boivin_breathing_2024,
	AUTHOR = {Boivin, Théo and Mathieu, Benoit and Porcher, Willy and Gillia, 
	Olivier},
	PUBLISHER = {IOP Publishing},
	DATE = {2024},
	DOI = {10.1149/1945-7111/ad14cf},
	ISSN = {1945-7111},
	JOURNALTITLE = {Journal of The Electrochemical Society},
	LANGID = {english},
	NUMBER = {1},
	PAGES = {010505},
	SHORTTITLE = {Breathing of a {{Silicon-Based Anode}}},
	TITLE = {Breathing of a {{Silicon-Based Anode}}: {{Mechanical Discrete 
	Approach Using DEM}}},
	URLDATE = {2024-01-09},
	VOLUME = {171},
}

@BOOK{cumberland_packing_1987,
	AUTHOR = {Cumberland, D. J. and Crawford, R. J.},
	PUBLISHER = {Elsevier Science Pub. Co. Inc., Amsterdam, NY},
	DATE = {1987},
	LANGID = {english},
	TITLE = {{{THE PACKING OF PARTICLES}}},
	URLDATE = {2024-10-14},
	URL = {https://openlibrary.org/books/OL15084426M/The_packing_of_particles},
}

@INPROCEEDINGS{bowers_introduction_2008,
	AUTHOR = {Bowers, Philip L.},
	BOOKTITLE = {Circle {{Packing}}: {{The Theory}} of {{Discrete Analytic 
	Functions}}},
	DATE = {2008},
	PAGES = {1--16},
	SHORTTITLE = {{{INTRODUCTION TO CIRCLE PACKING}}},
	TITLE = {Introduction to {{Circle Packing}}: {{A Review}}},
	URLDATE = {2024-11-13},
	URL = {https://api.semanticscholar.org/CorpusID:16173679},
}

@ARTICLE{hifi_literature_2009,
	AUTHOR = {Hifi, Mhand and M'Hallah, Rym},
	PUBLISHER = {Hindawi},
	DATE = {2009},
	DOI = {10.1155/2009/150624},
	ISSN = {1687-9147},
	JOURNALTITLE = {Advances in Operations Research},
	LANGID = {english},
	PAGES = {150624},
	SHORTTITLE = {A {{Literature Review}} on {{Circle}} and {{Sphere Packing 
	Problems}}},
	TITLE = {A {{Literature Review}} on {{Circle}} and {{Sphere Packing 
	Problems}}: {{Models}} and {{Methodologies}}},
	URLDATE = {2024-05-29},
	VOLUME = {2009},
}

@ARTICLE{lubachevsky_minimum_2009,
	AUTHOR = {Lubachevsky, Boris D. and Graham, Ronald L.},
	DATE = {2009},
	DOI = {10.1016/j.disc.2008.03.017},
	ISSN = {0012-365X},
	JOURNALTITLE = {Discrete Mathematics},
	NUMBER = {8},
	PAGES = {1947--1962},
	TITLE = {Minimum Perimeter Rectangles That Enclose Congruent 
Non-Overlapping 
	Circles},
	URLDATE = {2024-10-28},
	VOLUME = {309},
}

@ARTICLE{specht_high_2013,
	AUTHOR = {Specht, E.},
	DATE = {2013},
	DOI = {10.1016/j.cor.2012.05.011},
	ISSN = {0305-0548},
	JOURNALTITLE = {Computers \& Operations Research},
	NUMBER = {1},
	PAGES = {58--69},
	TITLE = {High Density Packings of Equal Circles in Rectangles with Variable 
	Aspect Ratio},
	URLDATE = {2024-10-28},
	VOLUME = {40},
}

@ARTICLE{markot_improved_2021,
	AUTHOR = {Markót, Mihály Csaba},
	DATE = {2021},
	DOI = {10.1007/s10898-021-01086-z},
	ISSN = {1573-2916},
	JOURNALTITLE = {Journal of Global Optimization},
	LANGID = {english},
	NUMBER = {3},
	PAGES = {773--803},
	TITLE = {Improved Interval Methods for Solving Circle Packing Problems in 
the 
	Unit Square},
	URLDATE = {2024-10-28},
	VOLUME = {81},
}

@ARTICLE{amore_circle_2023,
	AUTHOR = {Amore, Paolo},
	DATE = {2023},
	DOI = {10.1063/5.0140644},
	ISSN = {1070-6631},
	JOURNALTITLE = {Physics of Fluids},
	NUMBER = {2},
	PAGES = {027130},
	TITLE = {Circle Packing in Regular Polygons},
	URLDATE = {2024-10-28},
	VOLUME = {35},
}

@ARTICLE{gomadam_modeling_2005,
	AUTHOR = {Gomadam, Parthasarathy M. and Weidner, John W.},
	PUBLISHER = {IOP Publishing},
	DATE = {2005},
	DOI = {10.1149/1.2136087},
	ISSN = {1945-7111},
	JOURNALTITLE = {Journal of The Electrochemical Society},
	LANGID = {english},
	NUMBER = {1},
	PAGES = {A179},
	TITLE = {Modeling {{Volume Changes}} in {{Porous Electrodes}}},
	URLDATE = {2021-09-28},
	VOLUME = {153},
}

@MISC{chang_simple_2010,
	AUTHOR = {Chang, Hai-Chau and Wang, Lih-Chung},
	LOCATION = {Cornell University},
	DATE = {2010},
	DOI = {10.48550/arXiv.1009.4322},
	JOURNALTITLE = {arXiv e-prints},
	TITLE = {A {{Simple Proof}} of {{Thue}}'s {{Theorem}} on {{Circle 
Packing}}},
	URLDATE = {2024-05-29},
}

@MISC{gitlab_breathing,
	URL = {https://gitlab.com/TheoBoivin/articles/Breathing-Coefficient},
	TITLE = {{{GitLab}} -- Théo Boivin -- Articles -- {{Breathing 
	Coefficient}}},
	VERSION = {1.0.0},
}

\end{document}